\documentclass[12pt]{article}
\usepackage{fullpage, amssymb, amsmath, amsthm, mathabx}
\usepackage{color}
\usepackage{tikz}
\usetikzlibrary{shapes,arrows}
\usepackage{setspace}
\theoremstyle{plain}
\newtheorem{thm}{Theorem}

\newtheorem*{pigeonhole}{Pigeonhole Principle}
\newtheorem*{hall}{Hall's Marriage Theorem}
\newtheorem*{konig}{K\H{o}nig's Theorem}
\newtheorem*{vizing}{Vizing's Theorem}

\theoremstyle{definition}

\newtheorem*{atticus}{Player's move}
\newtheorem*{lazy_demon}{Lazy Demon}
\newtheorem*{contrary_demon}{Contrary Demon}
\newtheorem*{konig_demon}{K\H{o}nig's Demon}
\newtheorem*{vizing_demon}{Vizing's Demon}
\newtheorem*{winning}{Winning}

\theoremstyle{remark}

\newcommand{\set}[1]{\left\{ #1 \right\}}

\newcommand{\card}[1]{\left|#1\right|}

\newcommand{\parens}[1]{\left( #1 \right)}

\newcommand{\DefinedAs}{\mathrel{\mathop:}=}

\title{Playing cards with Vizing's demon\thanks{[Editorial Comment, added 2021] This paper was originally written in 2011 and updated in 2012. It was submitted to an expository journal but rejected, and never resubmitted. The second author posted a related article with some overlapping results ``A game generalizing Hall's theorem" on arXiv:1204.0139, and in 2014 published it in the journal \textit{Discrete Mathematics} (see \cite{rabern2012game}).}}
\author{Brian Rabern and Landon Rabern}
\date{}
\begin{document}
\maketitle
\begin{abstract}
\noindent We analyze a solitaire game in which a demon rearranges some cards after each move.  The graph edge coloring theorems of K\H{o}nig (1931) and Vizing (1964) follow from the winning strategies developed.
\end{abstract}

\section{Introduction}
For the uninitiated, the dense nature of mathematical language can act as an obscuring force.  With this essay we aim to bring two classical results of discrete mathematics into the light.  To this end we analyze winning strategies in a certain class of solitaire games.  The gains are non-standard proofs of the results of K\H{o}nig \cite{konig} and Vizing \cite{vizing}.  For the standard treatment of these results, see \cite{stiebitz}.  (For a dense and obscure version of the non-standard proofs presented here, see \cite{rabern2012game}.)  First, let's introduce the games.

\section{The Solitaire Game}

Before setting up a game a \textit{game number} $k$ is required. The game number corresponds to how many fixed \textit{stacks} that the game will employ and also determines other key features of the game.  We disallow the trivial case where there are no stacks (i.e. $k \geq 1$). Call a game with $k$ stacks a $k$-game.

The game also requires a \textit{demon}. To set up a $k$-game the demon first chooses a \textit{card number} $m\geq k$. The card number $m$ will correspond to how many distinctly numbered cards that will be involved (e.g. if $m = 3$, then the game will involve only 1-cards, 2-cards, and 3-cards). The demon, then, creates a deck with $k$ $1$-cards, $k$ $2$-cards, \ldots, and $k$ $m$-cards. With the game materials thus determined, the demon deals the stacks: the demon constructs $k$ nonempty stacks such that each stack contains at most one card of each type.  The total number of cards put in each stack is up to the demon. For each stack $i$ let $n_i$ be the number of cards the demon deals into that stack, then the sequence $(n_1, \ldots, n_k)$ represents the sizes of the stacks---call this sequence the game's \textit{stack profile}. The leftover cards in the deck form the \textit{reserve}. (The stacks and the reserve are displayed face up with all the cards showing.)

For example, consider a 3-game where the demon chooses a card number of 4 and a stack profile $(1,1,2)$. The demon will create 3 stacks of cards selected from 1-cards, 2-cards, 3-cards, or 4-cards such that the first stack contains one card, the second contains one card, and the the third stack contains two distinct cards (and the reserve contains the rest). Table \ref{table1}, then, represents the initial state of the game.

\begin{table}[ht]
\begin{center}

\begin{tabular}{c  c  c  c}
\hline
\multicolumn{1}{| r |}{Stack 1} & Stack 2 & \multicolumn{1}{| r |}{Stack 3} & \multicolumn{1}{|| r ||}{Reserve} \\
\hline \\
\colorbox{white}{\framebox{2}} & \framebox{2} & \framebox{2} & 
\colorbox{white}{\framebox{1}} \framebox{1} \framebox{1} \\
& & \framebox{4} & \framebox{3} \framebox{3} \\
& & \framebox{3} & \framebox{4} \framebox{4} \\

\end{tabular}
\caption{A 3-game set up.}
\end{center}
\label{table1}
\end{table}

With a game set up, the player repeatedly make moves of the following form.

\begin{atticus}
Pick some stack containing an $a$-card but no $b$-card and then swap
the $a$-card for a $b$-card from the reserve. In other words, the player may swap any card from a stack with some card in the reserve but only if this does not result in there being two cards of the same type in a stack.
\end{atticus}

The aim of the game is simply the following.

\begin{winning}
The player wins the game if at the start of a turn she can make a hand of $k$ differently numbered cards by picking one card from each stack.
\end{winning}

The demon, however, will complicate things by rearranging some cards after each of the player's moves.  Whether or not there is a winning strategy will depend on how the demon chooses to rearrange the cards.  Consider for example the following extreme demons.

\begin{lazy_demon}
After each turn, Lazy Demon does nothing.
\end{lazy_demon}

Can a player win against this demon?  Sure---all they'd have to do is go through the stacks in order, swapping an $i$-card into the $i$-th stack if there is not one there already. In the example given above the player need only swap a 1-card into Stack 1 to get the win.

\begin{table}[ht]
\begin{center}

\begin{tabular}{c  c  c  c}
\hline
\multicolumn{1}{| r |}{Stack 1} & Stack 2 & \multicolumn{1}{| r |}{Stack 3} & \multicolumn{1}{|| r ||}{Reserve} \\
\hline \\
\colorbox{green}{\framebox{1}} & \framebox{2} & \framebox{2} & \colorbox{green}{\framebox{2}} \framebox{1} \framebox{1} \\
& & \framebox{4} & \framebox{3} \framebox{3} \\
& & \framebox{3} & \framebox{4} \framebox{4} \\

\end{tabular}
\caption{A winning move.}
\end{center}
\label{table2}
\end{table}

\begin{contrary_demon}
After each turn, Contrary Demon undoes what what was just done, e.g. say the player swapped an $a$-card out of and a $b$-card into the $i$-th stack, then Contrary Demon swaps a $b$-card out of and an $a$-card into the $i$-th stack.
\end{contrary_demon}

When faced with this demon a player cannot ever change the stacks and so only wins if the demon gave her a winning position to start with.

\section{Some Principled Demons}

The extreme demons make the game boring but various intermediate demons provide for interesting games and strategy.  In fact, the graph-theoretic edge coloring theorems of K\H{o}nig and Vizing follow from the winning strategies developed.

For ease of illustration, let's say that Atticus is the solitaire master, who is playing the game. What we'll be concerned to demonstrate is that Atticus has a winning strategy when confronted with various principled demons.

\subsection{K\H{o}nig's Demon}

\begin{konig_demon}
Say Atticus swapped an $a$-card out of and a $b$-card into the $i$-th stack.
K\H{o}nig's Demon either does nothing or picks a stack other than the $i$-th containing a $b$-card but no $a$-card and swaps the $b$-card for an $a$-card from the reserve.
\end{konig_demon}

\begin{thm}
Atticus has a winning strategy against K\H{o}nig's Demon for any $k$-game.
\end{thm}

Atticus takes the following strategy.  Suppose the demon has set up a $k$-game with stack profile $(n_1, \ldots, n_k)$.  Atticus should consider the largest hand of differently numbered cards he can make by picking one card from each stack.  If he can make a hand of size $k$, he wins.  Otherwise, there is some stack, say the $i$-th, from which he is not picking a card and some number $b \leq m$ not appearing on any card in his hand.  Since $n_i \geq 1$, there is at least one card in the $i$-th stack, say an $a$-card.  Now Atticus can swap an $a$-card out of and a $b$-card into the $i$-th stack.  Atticus can now make a larger hand by picking the $b$-card from the $i$-th stack.  Since the hand uses a $b$-card from only the $i$-th stack, the demon swapping out a $b$-card in another stack cannot decrease the size of Atticus' hand.  Repeating this process, Atticus ends up with a hand of size $k$ and wins.

\subsection{Vizing's Demon}

\begin{vizing_demon}
Say Atticus' swapped an $a$-card out of and a $b$-card into the $i$-th stack.
Vizing's Demon either does nothing or picks a stack other than the $i$-th, containing a $b$-card but no $a$-card and swaps the $b$-card for an $a$-card from the deck, or picks a stack other than the $i$-th containing an $a$-card but no $b$-card and swaps the $a$-card for a $b$-card from the deck.
\end{vizing_demon}

\begin{thm}
Atticus has a winning strategy against Vizing's Demon for any $k$-game with stack profile $\parens{n_1, \ldots, n_k}$ where at most one of the $n_i$ is $1$.\footnote{The strategy Atticus will employ is based, in part, on Ehrenfeucht, Faber and Kierstead's  proof of Vizing's theorem \cite{Ehrenfeucht1984159} and Schrijver's proof of Vizing's theorem \cite{schrijver}.}
\end{thm}

Suppose that for any such $k$-game, Atticus has a strategy to get to a position where, for some nonempty subset $S$ of at most $k-1$ stacks, there is a choice of differently numbered cards (say with numbers $a_1$, \ldots, $a_s$), one from each stack in $S$, so that the numbers on these cards appear in none of the stacks outside $S$. We call such a position \emph{reducible} for the following reason. Let $T$ be the rest of the stacks, say $T$ contains $t$ stacks.  Now Atticus mentally removes all the $a_i$ from the reserve and plays the $t$-game on the stacks $T$.  Note that $T$ must still satisfy our requirement on the stack profile.  Also, a winning hand for $T$ can be put together with the cards chosen for $S$ to get a winning hand for the original game. By repeatedly applying this, Atticus gets down to a $1$-game which he easily wins and then constructs a winning hand for the original game as just mentioned.

We now describe Atticus' strategy for getting to a reducible position. Suppose there are at most $k-1$ numbers appearing on the cards in the stacks.  By our assumption on the stack profile $\parens{n_1, \ldots, n_k}$, there are at least $2(k-1) + 1 = 2k-1$ total cards in the stacks.  Another way to count the total number of cards is to add up the number of times each number appears on a card in the stacks.  Since there are at most $k-1$ different numbers on the cards, some number appears on $3$ or more cards on the stacks by the Pigeonhole Principle.

\begin{pigeonhole}
If N pigeons are put into M pigeonholes, then at least one pigeonhole must contain at least the round up of $\frac{N}{M}$ pigeons.
\end{pigeonhole}

Say it is the number $a$ that appears on at least three cards.  Since at most $k-1$ numbers appear on the stacks, Atticus can choose a number $b$ not appearing on any stack and swap an $a$-card out of and a $b$-card into one of the stacks.  Since no other stack contains a $b$-card, the demon can only pass his turn or pick a stack other than the $i$-th containing an $a$-card but no $b$-card and swap the $a$-card for a $b$-card from the deck.  Since there were at least $3$ $a$-cards in the stacks, after the demon's response there is still at least one $a$-card in the stacks and now there is a $b$-card in the stacks as well.  So Atticus has increased the number of numbers appearing on the cards in the stacks.  He can now repeat this procedure until at least $k$ numbers appear on the stacks.

We claim that any position with at least $k$ numbers appearing on the stacks is either reducible or Atticus can pick a winning hand.  Suppose we have such a position and let $A$ be $k$ numbers appearing on its stacks. Choose the smallest nonempty subset $B$ of $A$ so that the numbers in $B$ appear on at most $\card{B}$ stacks.  We can make this choice because $A$ itself is such a subset.  First, suppose $B$ has size one and let $s$ be the only stack containing a card with $B$'s element on it.  Then using $S \DefinedAs \set{s}$ shows that the position is reducible.  Therefore we may assume that $\card{B} \geq 2$. Choose $b$ in $B$ and remove it to get $B'$.  By our choice of $B$ as smallest, the numbers in $B'$ appear on at least $\card{B'} + 1 = \card{B}$ stacks.  Therefore the numbers in $B$ appear on exactly $\card{B}$ stacks, let $S$ be these stacks.  To reduce using $S$ we need to be able to pick a hand of differently numbered cards, one from each stack in $S$.  We can do this using Hall's Marriage Theorem from 1935 \cite{hall}.

\begin{hall}
Suppose a village has $n$ men and $n$ women.  For each woman there is a group of  men she would happily marry and any man would happily marry any woman.  Then the men and women in the village can be paired up in marriages so that everyone is happy if and only if for every group of women, there is a group of men of equal size in which each man has at least one woman in the group who will happily marry him.
\end{hall}

To applying this to our situation, let $S$ be the men and $B$ the women.  For each woman $b$ in $B$, let the men she will happily marry be the stacks in $S$ containing a card with $b$ on it.  We know that $S$ and $B$ are the same size and for any subset of $C$ of $B$ the numbers in $C$ appear on at least $\card{C}$ stacks.  That is, there are the same number of men as women and for any set of women there is a group of men of equal size in which each man has at least one woman in the group who will happily marry him.  Therefore, by Hall's Marriage Theorem, the men and women in the village can be paired up in marriages so that everyone is happy.  But that means precisely that we can pick a hand of differently numbered cards, one from each stack in $S$.

If $\card{S} = k$, then Atticus can pick a winning hand in the position.  Otherwise, the position is reducible.  This proves the claim. Therefore Atticus has a winning strategy against Vizing's Demon.  Figure \ref{fig:VizingFlow} shows this strategy.

\tikzstyle{decision} = [diamond, draw, fill=blue!20, 
    text width=4.5em, text badly centered, node distance=3cm, inner sep=4pt]
\tikzstyle{block} = [rectangle, draw, fill=blue!20, 
    text width=5em, text centered, rounded corners, minimum height=4em, inner sep=4pt]
\tikzstyle{line} = [draw, -latex']
\tikzstyle{cloud} = [draw, ellipse,fill=red!20, node distance=3cm,
    minimum height=2em]
    
\begin{figure}
\begin{center}
\begin{tikzpicture}[node distance = 4cm, auto]
    \node [block] (init) {initial position};
    \node [cloud, above of=init] (demon) {demon};
    \node [decision, right of=init, node distance = 4cm] (isWinning) {winning hand?};
    \node [decision, below of=isWinning, node distance = 4cm] (canReduce) {reducible?};
    \node [block, left of=canReduce] (reduce) {reduce to smaller game};
    \node [block, right of=canReduce] (increase) {increase count of numbers appearing on the stacks};
    \node [block, right of=isWinning] (win) {win};   

    \path [line,dashed] (demon) -- (init);
    \path [line] (init) -- (isWinning);
    \path [line] (isWinning) -- node [near start] {yes} (win);
    \path [line] (isWinning) -- node {no}(canReduce);

	 \path [line] (canReduce) -- node [near start] {yes} (reduce);
    \path [line] (canReduce) -- node {no}(increase);
    \path [line,dashed] (reduce) -- (init);
	 \path [line] (increase) -- (isWinning);
\end{tikzpicture}
\end{center}
\caption{Atticus' strategy for defeating Vizing's Demon}
\label{fig:VizingFlow}
\end{figure}
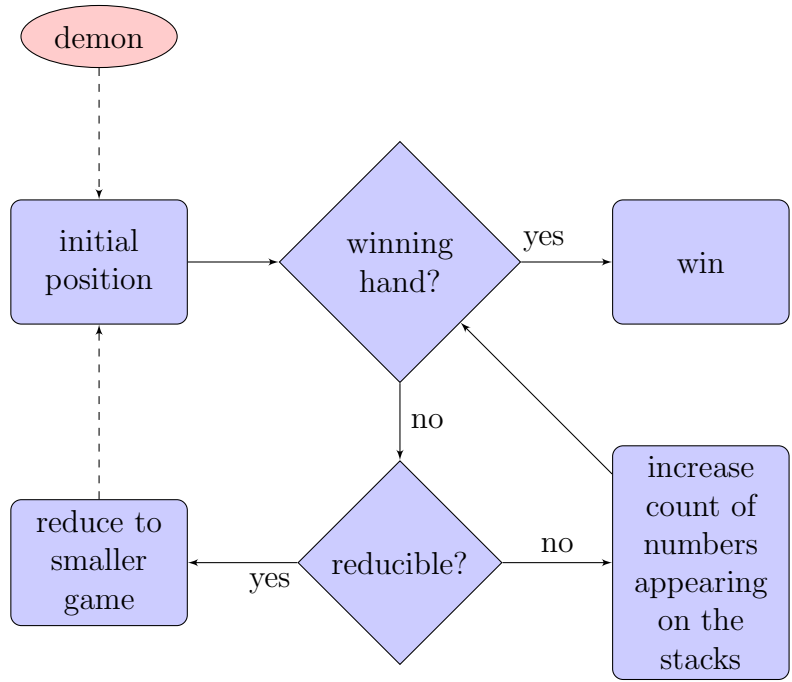

\section{Edge Coloring}\label{graphcolor}
A \emph{graph} consists of some dots and lines between them.  We call the dots \emph{vertices} and the lines \emph{edges}.  Also, we insist that no edge goes from a vertex to itself and at most one edge goes between any two vertices.  An \emph{edge coloring} of a graph is an assignment of colors to its edges so that no incident edges get the same color.  Two vertices are called \emph{neighbors} if there is an edge between them.  A graph is \emph{bipartite} if its vertices can be divided into two groups so that no two vertices in the same group are neighbors.  A \emph{cycle} in a graph is a sequence of different vertices such that each is a neighbor of the next in turn and the last is a neighbor of the first.  No bipartite graph contains a cycle with an odd number of vertices.

A classical theorem of K\H{o}nig from 1931 \cite{konig} can be written in the following form.

\begin{konig}
If every vertex in a bipartite graph has at most $k$ neighbors, then the graph has an edge coloring using $k$ colors.\footnote{It should be noted that a color can be used zero times, so if a graph can be edge colored using $k$ colors then it can be edge colored using $t$ colors for any $t$ bigger than $k$ as well.}
\end{konig}

Vizing's theorem from 1964 \cite{vizing} shows that only one more color suffices for general graphs.

\begin{vizing}
If every vertex in a graph has at most $k$ neighbors, then the graph has an edge coloring using $k + 1$ colors.
\end{vizing}

Both theorems can be reduced to a solitaire game as follows. Suppose there is a vertex $v$ in a graph $G$ with $k$ neighbors so that if we remove $v$ (and all edges incident to $v$) we can edge color the resulting graph using the colors $1, 2, \ldots, m$ where $m$ is at least $k$.  We wish to extend this edge coloring to an edge coloring of $G$ using only the colors $1, 2, \ldots, m$.  We will play a solitaire game with game number $k$ and card number $m$.  For each neighbor $x$ of $v$ the demon creates a stack $S_x$ with one card for each number in $1,2,\ldots, m$ that does not appear on an edge incident to $x$.  If $x$ has $d$ neighbors in $G$, then $S_x$ has $m + 1 - d$ cards since all of the edges incident to $x$ get different colors, except the edge to $v$ gets no color.  Now suppose $S_x$ contains an $a$-card but no $b$-card.  Then there is an edge incident to $x$ colored with $b$ but none colored with $a$.  Consider a path starting at $x$ and alternating between edges colored $b$ and edges colored $a$.  Since at most one edge incident to any vertex is colored with any given color, there is a unique such path that is longest.  If we swap colors $a$ and $b$ along this path we get another edge coloring of $G$ without $v$ using only the colors $1, 2, \ldots, m$.  Moreover, $S_x$ has been changed by swapping its $a$-card for a $b$-card from the reserve.  If the path does not end at a neighbor of $v$, then no other stack is changed (so the demon passed his turn).  So suppose the path does end at a neighbor $y$ of $v$.  Then $y$ is not $x$ because $x$ has no indident edge colored with $a$ and only one incident edge colored with $b$. If $G$ is bipartite, then the path must have an even number of edges for otherwise with the edges between $v$ and $x$ and $v$ and $y$ it would create a cycle with an odd number of vertices.  Therefore, since the path started with an edge colored $b$, it must end with an edge colored $a$.  Then after swapping along the path $S_y$ has been changed by swapping its $b$-card for a $a$-card from the reserve.  That is, swapping along the path gives precisely an Atticus move followed by a K\H{o}nig's Demon move.  If $G$ is not bipartite, then the final edge might be colored either $a$ or $b$ and hence swapping along the path corresponds precisely to an Atticus move followed by a Vizing's Demon move.  

Now suppose one of the theorems were false, then we could pick a graph $G$ with the least number of vertices for which it fails.  Since a graph with only one vertex has no edges, its edges can be colored using no colors.  Therefore $G$ has at least two vertices.  Suppose $G$ has a vertex $v$ with $k$ neighbors and no vertex with more than $k$ neighbors. By minimality of $G$, removing $v$ from $G$  gives a graph that can be colored using $m$ colors where $m$ is $k$ for K\H{o}nig's theorem and $m$ is $k+1$ for Vizing's theorem.  Then for any neighbor of $x$ of $v$, the stack $S_x$ has $m + 1 - d$ cards where $d$ is the number of neighbors $x$ has.  Since $d$ is at most $k$, this is at least $m + 1 - k$ cards.  That is, at least one card for K\H{o}nig's theorem and at least two cards for Vizing's theorem.  But then by our results on solitaire games, Atticus can make a hand of $k$ differently numbered cards by picking one card from each stack.  Coloring the edges incident to $v$ with the numbers on these cards gives a coloring of $G$ using $m$ colors.  This is a contradiction.  Therefore our assumption that one of the theorems was false must have been false.  That is, both theorems are true!


\bibliographystyle{amsplain}
\bibliography{demon}
\end{document}